\documentclass{amsart}
\usepackage{graphicx}
\usepackage{amscd}
\usepackage{amsmath}
\usepackage{amsfonts}
\usepackage{amssymb}
\newtheorem{theorem}{Theorem}
\theoremstyle{plain}

\newtheorem{corollary}{Corollary}

\newtheorem{definition}{Definition}

\newtheorem{lemma}{Lemma}

\newtheorem{remark}{Remark}

\numberwithin{equation}{section}

\begin{document}
\title[Sobolev Inequalities via Truncation]{Sobolev Inequalities: Symmetrization and Self Improvement via Truncation}
\author{Joaquim Martin$^{\ast}$}
\address{Department of Mathematics\\
Universidad Autonoma de Barcelona}
\email{jmartin@mat.uab.es}
\author{Mario Milman}
\address{Department of Mathematics\\
Florida Atlantic University}
\email{extrapol@bellsouth.net}
\urladdr{http://www.math.fau.edu/milman}
\author{Evgeniy Pustylnik}
\address{Department of Mathematics\\
Technion, Israel Institute of Technology}
\email{evg@techunix.technion.ac.il}
\thanks{2000 Mathematics Subject Classification Primary: 46E30, 26D10.}
\thanks{$^{\ast}$ Supported in part by MTM2004-02299 and by CURE 2005SGR00556}
\thanks{This paper is in final form and no version of it will be submitted for
publication elsewhere.}

\begin{abstract}
We develop a new method to obtain symmetrization inequalities of Sobolev type.
Our approach leads to new inequalities and considerable simplification in the
theory of embeddings of Sobolev spaces based on rearrangement invariant spaces.
\end{abstract}\maketitle

\section{Introduction}

A well known basic principle in the theory of Sobolev embeddings, due to
Maz'ya, and Federer and Fleming (cf. \cite{Ma} and the references therein), is
the equivalence\footnote{In fact, the equivalence is sharp all the way down to
the constants.} between the isoperimetric inequality and the
Gagliardo-Nirenberg inequality
\begin{equation}
\left\|  f\right\|  _{L^{n/(n-1)}}\leq c_{n}\left\|  \nabla f\right\|
_{L^{1}},\forall f\in C_{0}^{\infty}(\mathbb{R}^{n}). \label{uno}%
\end{equation}

A second, somewhat less well know principle, which is often rediscovered in
the literature\footnote{See \cite{Ha}.}, and is also apparently due to Maz'ya
\cite{Ma}, states that, roughly speaking, under rather general circumstances a
weak type Sobolev inequality implies a strong type Sobolev inequality. We
refer to \cite{BaCLS}, \cite{Tar} and \cite{Ha}. In particular, the first two
quoted papers show how weak $L^{p}$ Sobolev inequalities self improve by
truncation to $L(p,q)$ inequalities, while \cite{Ha} provides a nice survey
and a unified treatment of the cases $p=1$ and $1<p<n,$ of the Sobolev embedding.

It is also known that Sobolev inequalities have an in-built *reiteration*
property which is due to a combination of the chain rule and H\"{o}lder's
inequalities. For example, since for any $\alpha>1$ we have $\left|
\nabla\left|  f\right|  ^{\alpha}\right|  =\alpha\left|  f\right|  ^{\alpha
-1}\left|  \nabla f\right|  ,$ it follows that if we pick $p\in(1,n),$ and let
$q=\frac{np}{n-p},\alpha=\frac{(n-1)p}{n-p}=\frac{n-1}{n}q,$ we have
$(\alpha-1)p^{\prime}=q,$ $q(\frac{n-1}{n}-\frac{1}{p^{\prime}})=1,$ and
$\left\|  f\right\|  _{L^{q}}^{q}=\left\|  \left|  f\right|  ^{\alpha
}\right\|  _{L^{n/(n-1)}}^{n/(n-1)}.$ Therefore, from (\ref{uno}) we thus have
that, for $f\in C_{0}^{\infty}(\mathbb{R}^{n})$,
\begin{align*}
\left\|  f\right\|  _{L^{q}}^{q(n-1)/n}  &  =\left\|  \left|  f\right|
^{\alpha}\right\|  _{L^{n/(n-1)}}\leq c_{n}\left\|  \alpha\left|  f\right|
^{\alpha-1}\left|  \nabla f\right|  \right\|  _{L^{1}}\\
&  \leq c_{n}\alpha\left\|  f\right\|  _{L^{q}}^{q/p^{\prime}}\left\|  \nabla
f\right\|  _{L^{p}},
\end{align*}
which immediately yields the classical Sobolev inequality.

It follows from the discussion above that, roughly speaking, ``all'' $L^{p}$
Sobolev inequalities follow the Gagliardo-Nirenberg inequality (\ref{uno}) or,
equivalently, from the isoperimetric inequality. But one can go further.
Talenti \cite{Ta}, using the isoperimetric inequality and the co-area formula,
obtained a powerful rearrangement inequality\footnote{For a related inequality
see also \cite{Ma}, Lemma 2.3.3.}, which is very close to the
P\'{o}lya-Szeg\"{o} principle (cf. (\ref{PSz}) below)
\begin{equation}
s^{1-1/n}\left(  -f^{\ast}\right)  ^{^{\prime}}(s)\leq c_{n}\frac{\partial
}{\partial s}\int_{\left\{  \left|  f\right|  >f^{\ast}(s)\right\}  }\left|
\nabla f(x)\right|  dx,\label{dos}%
\end{equation}
where $f^{\ast}$ denotes the non-increasing rearrangement of $f.$

In particular, Talenti's inequality can be used to prove Sobolev inequalities
in the setting of rearrangement invariant spaces (cf. \cite{Ta}, \cite{EKP}),
where, in principle, the chain rule argument is not available. Moreover, given
the precise information about the constant $c_{n}$ in (\ref{dos}), Talenti's
inequality allows one to obtain best possible constants for the classical
Sobolev inequalities (cf. \cite{Ta}).

A somewhat different rearrangement inequality\footnote{A slightly different
but equivalent inequality had been obtained earlier in \cite{Ko}.} was used in
\cite{BMR} to study the borderline case $p=n$,
\begin{equation}
f^{\ast\ast}(t)-f^{\ast}(t)\leq c_{n}t^{1/n}\left|  \nabla f\right|
^{\ast\ast}(t),f\in C_{0}^{\infty}(\mathbb{R}^{n}),t>0,\label{tres}%
\end{equation}
where $f^{\ast\ast}(t)=\frac{1}{t}\int_{0}^{t}f^{\ast}(s)ds.$ The proofs of
(\ref{tres}) in \cite{BMR} and in \cite{Ko} use the symmetrization
principle\footnote{$f^{\circ}(x)=f^{\ast}(\gamma_{n}\left|  x\right|  ^{n}),$
is the symmetric decreasing rearrangmeent of $f,$ $\gamma_{n}$ is the measure
of the unit ball in $\mathbb{R}^{n}.$} of P\'{o}lya-Szeg\"{o},
\begin{equation}
\left|  \nabla f^{\circ}\right|  ^{\ast\ast}(t)\leq\left|  \nabla f\right|
^{\ast\ast}(t),f\in C_{0}^{\infty}(\mathbb{R}^{n}).\label{PSz}%
\end{equation}
The inequality (\ref{tres}) is further extended in \cite{MM} using both
Talenti's inequality (\ref{dos}) and the isoperimetric inequality.

The sharpest form of the classical Sobolev inequalities, including the
critical exponent $p=n,$ follow from (\ref{tres}), namely, for $1<p\leq
n,1\leq q\leq\infty,$ we have
\begin{equation}
\left\{  \int_{0}^{\infty}[(f^{\ast\ast}(t)-f^{\ast}(t))t^{1/p-1/n}]^{q}%
\frac{dt}{t}\right\}  ^{1/q}\leq c_{n,p}\left\{  \int_{0}^{\infty}[\left|
\nabla f\right|  ^{\ast}(t)t^{1/p}]^{q}\frac{dt}{t}\right\}  ^{1/q}.
\label{gardel}%
\end{equation}
It turns out, however, that the important case $p=1,$ which is also valid,
requires a separate argument since $c_{n,p}$ blows up as p tends to $1$.
Indeed (\ref{gardel}) for $p=1$ is the sharp form of the Gagliardo-Nirenberg
inequality due to Poornima \cite{Po} (cf. (\ref{teoA5}) below).

Symmetrization inequalities imply Sobolev inequalites in the setting of
rearrangement invariant spaces. Indeed, from (\ref{tres}) we obtain: for any
r.i. space $X$ with upper Boyd\footnote{The restriction on the Boyd indices is
only required to guarantee that the inequality $\left\|  g^{\ast\ast}\right\|
_{X}\leq c_{X}\left\|  g\right\|  _{X},$ holds for all $g\in X.$} index
$\beta_{X}<1$, we have (cf. \cite{MP})
\[
\left\|  t^{-1/n}(f^{\ast\ast}(t)-f^{\ast}(t))\right\|  _{X}\leq c\left\|
\nabla f\right\|  _{X},f\in C_{0}^{\infty}(\mathbb{R}^{n}),
\]
where $c=c(n,X).$ Moreover, the inequality is sharp (cf. Section \ref{evg}
below): if $Y$ is any r.i. space then the validity of%

\[
\left\|  f\right\|  _{Y}\leq c\left\|  \nabla f\right\|  _{X},f\in
C_{0}^{\infty}(\mathbb{R}^{n})
\]
implies that
\[
\left\|  f\right\|  _{Y}\leq\left\|  t^{-1/n}(f^{\ast\ast}(t)-f^{\ast
}(t))\right\|  _{X}.
\]

Note that for $X=L^{p}$ the condition $\beta_{X}<1$ translates into $p>1.$ The
fact that spaces near $L^{1}$ cannot be treated using (\ref{tres}), and the
previous discussion showing the central role of the Gagliardo-Nirenberg
inequalities [cf. (\ref{uno}) and (\ref{teoA5})], suggested that there could
be another more powerful underlying rearrangement inequality that would allow
for a unified treatment.

The purpose of this paper is to show that truncation can be actually used as a
method to obtain symmetrization inequalities. In other words rather than show
that a Sobolev inequality implies other Sobolev inequalities one case at a
time, we prove that from a Sobolev inequality we can obtain a symmetrization
inequality that ``implies all the Sobolev inequalities''.

Our analysis leads indeed to new symmetrization inequalities that allow for a
unified treatment of the Sobolev inequalites at both end points in the setting
of r.i. spaces. Remarkably, our approach also provides a considerable
simplification to the methods used to prove the classical symmetrization
inequalities discussed above. This is important for the application of our
methods to generalized settings like metric spaces (cf. \cite{KM}), fractional
derivatives (cf. \cite{MM1}), capacities, etc, which we hope to treat elsewhere.

The following is our main result. We could call it a ``symmetrization by
truncation principle'', and it is part of a family of similar inequalities, we
consider here the most important case, namely the end point $p=1$ (cf. Section
\ref{markao} below).

\begin{theorem}
\label{teoA}The following statements are equivalent

\begin{enumerate}
\item [(i)]%
\begin{equation}
W_{0}^{1,1}(\mathbb{R}^{n})\subset L^{n/(n-1),\infty}(\mathbb{R}^{n}).
\label{teoA1}%
\end{equation}

\item[(ii)]
\begin{equation}
\int_{0}^{t}s^{-\frac{1}{n}}[f^{\ast\ast}(s)-f^{\ast}(s)]ds\leq cn\int_{0}%
^{t}\left|  \nabla f\right|  ^{\ast}(s)ds,\ f\in C_{0}^{\infty}(\mathbb{R}%
^{n}). \label{teoA2}%
\end{equation}

\item[(iii)] For any rearrangement invariant space $X$ with lower Boyd
index\footnote{For $X=L^{p},\alpha_{L^{p}}=1/p>0$ translates into $p<\infty.$}
$\alpha_{X}>0 $ we have
\begin{equation}
\left\|  s^{-1/n}(f^{\ast\ast}(s)-f^{\ast}(s))\right\|  _{X}\preceq\left\|
\left|  \nabla f\right|  \right\|  _{X},f\in C_{0}^{\infty}(\mathbb{R}^{n}).
\label{final}%
\end{equation}

\item[(iv)]
\begin{equation}
W_{0}^{1,1}(\mathbb{R}^{n})\subset L^{n/(n-1),1}(\mathbb{R}^{n}).
\label{teoA5}%
\end{equation}
\end{enumerate}
\end{theorem}

To understand how Theorem \ref{teoA} represents an improvement over the known
results, we note that the implication (\ref{teoA1})$\Rightarrow$(\ref{teoA5})
is the self improvement that follows by the usual method of truncation (cf.
\cite{BaCLS}, \cite{Ta}, \cite{Ha}). On the other hand, by ``symmetrization by
truncation'' we obtain the new rearrangement inequality (\ref{teoA2}) which
readily gives (\ref{final}), and thus we have obtained the most general form
of the Sobolev inequalities in the context of r.i. spaces. Moreover, in the
process we have eliminated the restriction on the upper Boyd indices of
\cite{MP} and we are able to treat spaces near $L^{1}$ in a unified manner. In
particular, we note that Theorem \ref{teoA}, and the discussion preceding it,
shows that the symmetrization inequality\footnote{We shall also refer
sometimes to inequalities involving the quantity $f^{\ast\ast}(t)-f^{\ast}(t)$
as ``oscillation inequalities''.} (\ref{teoA2}) is equivalent to the
isoperimetric inequality.

Furthermore, since we believe that it is methodologically important for
further extensions, and in order to clarify the role of the assumptions that
intervene in the proof of the basic inequalities, in Section \ref{marka} below
we provide a simple direct proof that all the main rearrangement inequalities
discussed here, namely (\ref{dos}), (\ref{tres}) (\ref{teoA2}), (\ref{PSz}%
)\footnote{Actually the version we prove of (\ref{PSz}) is slightly weaker in
as much as the constant $n$ appears on the right hand side of the inequality.}
follow directly from the straightforward weak type Sobolev inequality
(\ref{teoA1}) via truncation.

A complete discussion concerning Sobolev embeddings in the setting of r.i.
spaces is then given in the final Section \ref{evg}. Our approach treats all
cases in a unified manner with optimal conditions, the optimal spaces are
explicitly constructed and, moreover, we give a unified treatment of all the
borderline cases as well (the reader should compare our approach with the ones
that are currently available in the literature: cf. \cite{EKP}, \cite{MP},
\cite{KP1}, \cite{KP2}, and the references quoted therein). We also show how
our methods provide a considerable simplification to recent results on the
compactness of Sobolev embeddings (cf. \cite{KP3} and \cite{P1}).

We stress that in this paper we have not attempted to prove the most general
results, but rather we aim to illustrate the power of our methods. In
particular, in order not to obscure the simplicity of the arguments we work
for the most part on $\mathbb{R}^{n},$ and we formulate our results as
inequalities. This is justified since the extensions to regular domains can be
obtained using well known techniques, while more sophisticated extensions
would require a separate treatment.

As usual, the symbol $f\simeq g$ will indicate the existence of a universal
constant $c>0$ (independent of all parameters involved) so that $(1/c)f\leq
g\leq c\,f$, while the symbol $f\preceq g$ means that $f\leq c\,g,$ and
$f\succeq g$ means that $f\geq c\,g.$

\section{Symmetrization Inequalities by Truncation\label{marka}}

The purpose of this section is to show that all the symmetrization
inequalities discussed in the Introduction follow from the Sobolev embedding
\begin{equation}
W_{0}^{1,1}(\mathbb{R}^{n})\subset L^{n/(n-1),\infty}(\mathbb{R}^{n}),
\label{wtyp}%
\end{equation}
by truncation.

Since it will be important for us to keep track of the constants of the
embedding (\ref{wtyp}), and in order to provide a self contained presentation,
we present a proof of (\ref{wtyp}) following \cite{Ha}, who in turn credits
Santalo for the method of proof.

\begin{lemma}
Let $f\in W_{0}^{1,1}(\mathbb{R}^{n}),$ then
\[
\sup_{t>0}t\left|  \left\{  x\in\mathbb{R}^{n}:\left|  f(x)\right|
>t\right\}  \right|  ^{\frac{n-1}{n}}\leq\frac{1}{\gamma_{n}^{1/n}}%
\int_{\mathbb{R}^{n}}\left|  \nabla f(x)\right|  dx,
\]
where $\gamma_{n}=$measure of the unit ball in $\mathbb{R}^{n}.$
\end{lemma}

\begin{proof}
Let $f\in C_{0}^{\infty}(\mathbb{R}^{n}),$ then as it is well known (see
\cite[Page 125]{St}) we have the representation
\[
f(x)=\frac{1}{n\gamma_{n}}\sum_{j=1}^{n}\int_{\mathbb{R}^{n}}\frac{\partial
f}{\partial x_{j}}(x-y)\frac{y_{j}}{\left|  y\right|  ^{n}}dy.
\]
Thus,
\[
\left|  f(x)\right|  \leq\frac{1}{n\gamma_{n}}\int_{\mathbb{R}^{n}}\left|
\nabla f(y)\right|  \frac{1}{\left|  x-y\right|  ^{n-1}}dy.
\]
Let $H=\left\{  x:\left|  f(x)\right|  >t\right\}  ,$ then, combining the
previous inequality with Chebyshev's inequality and Fubini, we find that
\[
t\left|  H\right|  \leq\int_{H}\left|  f(x)\right|  dx\leq\frac{1}{n\gamma
_{n}}\int_{\mathbb{R}^{n}}\left|  \nabla f(y)\right|  \int_{H}\frac
{dx}{\left|  x-y\right|  ^{n-1}}dy.
\]
For a fixed $y$ let $B=B(y,r)$ be a ball such that such $\left|  B\right|
=\left|  H\right|  .$ Then by symmetrization
\[
\int_{H}\frac{dx}{\left|  x-y\right|  ^{n-1}}\leq\int_{B}\frac{dx}{\left|
x-y\right|  ^{n-1}}=n\gamma_{n}r=n\gamma_{n}^{1-1/n}\left|  H\right|  ^{1/n}.
\]
Summarizing
\[
t\left|  H\right|  \leq\frac{\left|  H\right|  ^{1/n}}{\gamma_{n}^{1/n}}%
\int_{\mathbb{R}^{n}}\left|  \nabla f(y)\right|  dy.
\]
\end{proof}

\subsection{Talenti's inequality\label{seta}}

Our starting point is the weak type inequality
\begin{equation}
\sup_{t>0}t\left|  \left\{  x\in\mathbb{R}^{n}:\left|  f(x)\right|
>t\right\}  \right|  ^{\frac{n-1}{n}}\leq\gamma_{n}^{-1/n}\int_{\mathbb{R}%
^{n}}\left|  \nabla f(x)\right|  dx,f\in W_{0}^{1,1}(\mathbb{R}^{n}).
\label{dosdos}%
\end{equation}

Let $0<t_{1}<t_{2}<\infty,$ the truncations of $f$ are defined by
\[
f_{t_{1}}^{t_{2}}(x)=\left\{
\begin{array}
[c]{ll}%
t_{2}-t_{1} & \text{if }\left|  f(x)\right|  >t_{2},\\
\left|  f(x)\right|  -t_{1} & \text{if }t_{1}<\left|  f(x)\right|  \leq
t_{2},\\
0 & \text{if }\left|  f(x)\right|  \leq t_{1}.
\end{array}
\right.
\]
Observe that if $f\in W_{0}^{1,1}(\mathbb{R}^{n})$ then $f_{t_{1}}^{t_{2}}\in
W_{0}^{1,1}(\mathbb{R}^{n}),$ therefore replacing $f$ by $f_{t_{1}}^{t_{2}}$
in (\ref{dosdos}) we obtain
\[
\sup_{t>0}t\left|  \left\{  x\in\mathbb{R}^{n}:\left|  f_{t_{1}}^{t_{2}%
}(x)\right|  >t\right\}  \right|  ^{\frac{n-1}{n}}\leq\gamma_{n}^{-1/n}%
\int_{\mathbb{R}^{n}}\left|  \nabla f_{t_{1}}^{t_{2}}(x)\right|  dx.
\]
We obviously have
\[
\sup_{t>0}t\left|  \left\{  x\in\mathbb{R}^{n}:\left|  f_{t_{1}}^{t_{2}%
}(x)\right|  >t\right\}  \right|  ^{\frac{n-1}{n}}\geq\left(  t_{2}%
-t_{1}\right)  \left|  \left\{  x\in\mathbb{R}^{n}:\left|  f(x)\right|  \geq
t_{2}\right\}  \right|  ^{\frac{n-1}{n}},
\]
and
\[
\left|  \nabla f_{t_{1}}^{t_{2}}\right|  =\left|  \nabla f\right|
\chi_{\left\{  t_{1}<\left|  f\right|  \leq t_{2}\right\}  }.
\]
Therefore,
\[
\left(  t_{2}-t_{1}\right)  \left|  \left\{  x\in\mathbb{R}^{n}:\left|
f(x)\right|  \geq t_{2}\right\}  \right|  ^{1-1/n}\leq\gamma_{n}^{-1/n}%
\int_{\left\{  t_{1}<\left|  f\right|  \leq t_{2}\right\}  }\left|  \nabla
f(x)\right|  dx.
\]
Let $0\leq a<b,$ and consider $t_{1}=f^{\ast}(b),$ $t_{2}=f^{\ast}(a).$ Then
\begin{align}
\left(  f^{\ast}(a)-f^{\ast}(b\right)  )a^{1-1/n}  &  \leq\left(  f^{\ast
}(a)-f^{\ast}(b\right)  )\left|  \left\{  x\in\mathbb{R}^{n}:\left|
f(x)\right|  \geq f^{\ast}(a)\right\}  \right|  ^{1-1/n}\label{truncacion}\\
&  \leq\gamma_{n}^{-1/n}\int_{\left\{  f^{\ast}(b)<\left|  f\right|  \leq
f^{\ast}(a)\right\}  }\left|  \nabla f(x)\right|  dx\nonumber\\
&  \leq\gamma_{n}^{-1/n}\int_{0}^{b-a}\left|  \nabla f\right|  ^{\ast
}(s)ds,\nonumber
\end{align}
whence $f^{\ast}$ is locally absolutely continuous.

Let $s>0$ and $h>0;$ the previous considerations with $t_{1}=f^{\ast}(s+h)$
and $t_{2}=f^{\ast}(s)$ yield
\[
\left(  f^{\ast}(s)-f^{\ast}(s+h)\right)  s^{1-1/n}\leq\gamma_{n}^{-1/n}%
\int_{\left\{  f^{\ast}(s+h)<\left|  f\right|  \leq f^{\ast}(s)\right\}
}\left|  \nabla f(x)\right|  dx.
\]
Thus,
\[
\frac{\left(  f^{\ast}(s)-f^{\ast}(s+h)\right)  }{h}s^{1-1/n}\leq\frac
{\gamma_{n}^{-1/n}}{h}\int_{\left\{  f^{\ast}(s+h)<\left|  f\right|  \leq
f^{\ast}(s)\right\}  }\left|  \nabla f(x)\right|  dx.
\]
Letting $h\rightarrow0$ we obtain (\ref{dos}).

\subsection{The Oscillation Inequality}

We now prove the oscillation inequality (\ref{tres}). We will integrate by
parts, so let us note first that using (\ref{truncacion}) we have, for
$0<s<t,$%
\begin{equation}
s\left(  f^{\ast}(s)-f^{\ast}(t\right)  )\leq\gamma_{n}^{-1/n}s^{1/n}\int
_{0}^{t-s}\left|  \nabla f\right|  ^{\ast}(s)ds. \label{boca}%
\end{equation}
Now,%

\begin{align}
f^{\ast\ast}(t)-f^{\ast}(t)  &  =\frac{1}{t}\int_{0}^{t}\left(  f^{\ast
}(s)-f^{\ast}(t)\right)  ds\label{cua}\\
&  =\frac{1}{t}\left\{  \left[  s\left(  f^{\ast}(s)-f^{\ast}(t)\right)
\right]  _{0}^{t}+\int_{0}^{t}s\left(  -f^{\ast}\right)  ^{^{\prime}%
}(s)ds\right\} \nonumber\\
&  =\frac{1}{t}\int_{0}^{t}s\left(  -f^{\ast}\right)  ^{^{\prime}%
}(s)ds,\nonumber
\end{align}
where the integrated term $\left[  s\left(  f^{\ast}(s)-f^{\ast}(t)\right)
\right]  _{0}^{t}$ vanishes on account of (\ref{boca}).

Now, starting from (\ref{cua}) we readily get
\begin{align*}
f^{\ast\ast}(t)-f^{\ast}(t)  &  =\frac{1}{t}\int_{0}^{t}s\left(  -f^{\ast
}\right)  ^{^{\prime}}(s)ds=\frac{1}{t}\int_{0}^{t}s^{1/n}s^{1-1/n}\left(
-f^{\ast}\right)  ^{^{\prime}}(s)ds\\
&  \leq\frac{t^{1/n}}{t}\int_{0}^{t}s^{1-1/n}\left(  -f^{\ast}\right)
^{^{\prime}}(s)ds\\
&  \leq\gamma_{n}^{-1/n}\frac{t^{1/n}}{t}\int_{0}^{t}\left(  \frac{\partial
}{\partial s}\int_{\left\{  \left|  f\right|  >f^{\ast}(s)\right\}  }\left|
\nabla f(x)\right|  dx\right)  ds\\
&  \leq\gamma_{n}^{-1/n}t^{1/n}\left|  \nabla f\right|  ^{\ast\ast}(t),
\end{align*}
where in the third step we used (\ref{dos}).

\begin{remark}
Since it will be useful below we observe that in an intermediate step of the
previous derivation we implicitly obtained the inequality
\begin{equation}
\int_{0}^{t}s^{1-1/n}\left(  -f^{\ast}\right)  ^{^{\prime}}(s)ds\leq\gamma
_{n}^{-1/n}\int_{0}^{t}\left|  \nabla f\right|  ^{\ast}(s)ds.
\label{intermedia}%
\end{equation}
\end{remark}

\subsection{Integrated Oscillation Inequality\label{markao}}

We prove (\ref{teoA2}). Starting from (\ref{cua}) and integrating by parts we
have
\begin{align*}
\int_{0}^{t}s^{-\frac{1}{n}}[f^{\ast\ast}(s)-f^{\ast}(s)]ds  &  =\int_{0}%
^{t}s^{-1-\frac{1}{n}}\int_{0}^{s}u\left(  -f^{\ast}\right)  ^{^{\prime}%
}(u)du\text{ }ds\\
&  =-n\int_{0}^{t}\int_{0}^{s}u\left(  -f^{\ast}\right)  ^{^{\prime}%
}(u)du\text{ }ds^{-1/n}\\
&  =\left.  -ns^{-1/n}\int_{0}^{s}u\left(  -f^{\ast}\right)  ^{^{\prime}%
}(u)du\text{ }\right|  _{0}^{t}+\int_{0}^{t}s^{1-1/n}\left(  -f^{\ast}\right)
^{^{\prime}}(s)ds.
\end{align*}

Since by (\ref{boca}) and (\ref{cua}) it follows that
\[
s^{-1/n}\int_{0}^{s}u\left(  -f^{\ast}\right)  ^{^{\prime}}(u)du=s^{1-1/n}%
\left(  f^{\ast\ast}(s)-f^{\ast}(s)\right)  \preceq\int_{0}^{s}\left|  \nabla
f\right|  ^{\ast}(s)ds,
\]
the integrated term vanishes at $t=0.$ Consequently, in view of (\ref{tres})
and (\ref{intermedia}), we can continue our estimates with
\begin{align*}
\int_{0}^{t}s^{-\frac{1}{n}}[f^{\ast\ast}(s)-f^{\ast}(s)]ds  &  =-nt^{1-1/n}%
\left(  f^{\ast\ast}(t)-f^{\ast}(t)\right)  +\gamma_{n}^{-1/n}n\int_{0}%
^{t}\left|  \nabla f\right|  ^{\ast}(s)ds\\
&  \leq\gamma_{n}^{-1/n}n\int_{0}^{t}\left|  \nabla f\right|  ^{\ast}(s)ds,
\end{align*}
as we wished to show.

\begin{remark}
Using a standard limiting argument we may extend the validity of (\ref{tres})
and (\ref{teoA2}) from functions in $C_{0}^{\infty}(\mathbb{R}^{n})$ to all
functions in $W_{0}^{1,1}(\mathbb{R}^{n}).$ For example, suppose that
(\ref{teoA2}) holds for functions in $C_{0}^{\infty}(\mathbb{R}^{n})$. Then,
given $f\in W_{0}^{1,1}(\mathbb{R}^{n})$ select $f_{k}\in C_{0}^{\infty
}(\mathbb{R}^{n})$ such that
\[
f_{k}(x)\rightarrow f(x)\text{ a.e. and }f_{k}\rightarrow f\text{ in }%
W_{0}^{1,1}(\mathbb{R}^{n}).
\]
Since $f_{k}^{\ast}(t)\rightarrow f^{\ast}(t)$ a.e. we can use Fatou's lemma
\begin{align*}
\int_{0}^{t}s^{-\frac{1}{n}}[f^{\ast\ast}(s)-f^{\ast}(s)]ds  &  \leq\lim
\int_{0}^{t}s^{-\frac{1}{n}}[f_{k}^{\ast\ast}(s)-f_{k}^{\ast}(s)]ds\preceq
\lim\int_{0}^{t}\left|  \nabla f_{k}\right|  ^{\ast}(s)ds\\
&  =\lim\int_{0}^{t}\left|  \nabla(f_{k}+f-f)\right|  ^{\ast}(s)ds\\
&  \leq\lim\int_{0}^{t}\left|  \nabla(f_{k}-f)\right|  ^{\ast}(s)ds+\int
_{0}^{t}\left|  \nabla f\right|  ^{\ast}(s)ds\\
&  \leq\lim_{n}\left\|  \left|  \nabla(f_{k}-f)\right|  \right\|  _{L^{1}%
}+\int_{0}^{t}\left|  \nabla f\right|  ^{\ast}(s)ds\\
&  =\int_{0}^{t}\left|  \nabla f\right|  ^{\ast}(s)ds,
\end{align*}
as we wished to prove. The extension of (\ref{tres}) is proved similarly.
\end{remark}

\subsection{An elementary proof of the P\'{o}lya-Szeg\"{o} principle}

We will actually prove a slightly weaker form of the P\'{o}lya-Szeg\"{o}
principle, namely
\[
\left|  \nabla f^{\circ}\right|  ^{\ast\ast}(s)\leq n\left|  \nabla f\right|
^{\ast\ast}(s).
\]

Our starting point is Talenti's inequality (cf. Section \ref{seta} above): if
$f\in W_{0}^{1,1}(\mathbb{R}^{n})$ then
\[
s^{1-1/n}\left(  -f^{\ast}\right)  ^{^{\prime}}(s)\leq\gamma_{n}^{-1/n}%
\frac{\partial}{\partial s}\int_{\left\{  \left|  f\right|  >f^{\ast
}(s)\right\}  }\left|  \nabla f(x)\right|  dx.
\]

We claim that if $\Phi$ is a positive Young's function, then
\begin{equation}
\Phi\left(  n\gamma_{n}^{1/n}s^{1-1/n}\left(  -f^{\ast}\right)  ^{^{\prime}%
}(s)\right)  \leq\frac{\partial}{\partial s}\int_{\left\{  \left|  f\right|
>f^{\ast}(s)\right\}  }\Phi(n\left|  \nabla f(x)\right|  )dx.\label{pol}%
\end{equation}
Assuming momentarily the validity of (\ref{pol}) we get
\[
\int_{0}^{\infty}\Phi\left(  n\gamma_{n}^{1/n}s^{1-1/n}\left(  -f^{\ast
}\right)  ^{^{\prime}}(s)\right)  ds\leq\int_{R^{n}}\Phi(n\left|  \nabla
f(x)\right|  )dx,
\]
and since,
\[
\int_{0}^{\infty}\Phi\left(  n\gamma_{n}^{1/n}s^{1-1/n}\left(  -f^{\ast
}\right)  ^{^{\prime}}(s)\right)  ds=\int_{\mathbb{R}^{n}}\Phi(\left|  \nabla
f^{\circ}(x)\right|  )dx
\]
it follows that for all Young functions $\Phi$ we have
\[
\int_{\mathbb{R}^{n}}\Phi(\left|  \nabla f^{\circ}(x)\right|  )dx\leq
\int_{\mathbb{R}^{n}}\Phi(n\left|  \nabla f(x)\right|  )dx.
\]
The last inequality implies, by a well known result of
Hardy-Littlewood-P\'{o}lya (cf. \cite[Page 88]{BS}),
\[
\int_{0}^{t}\left|  \nabla f^{\circ}\right|  ^{\ast}(s)ds\leq n\int_{0}%
^{t}\left|  \nabla f\right|  ^{\ast}(s)ds,
\]
as we wished to show.

It remains to prove (\ref{pol}). Here we follow Talenti's argument (it is
important for our purposes to note that at this point in the argument we are
not using the isoperimetric inequality or the co-area formula)\textsl{. }Let
$s>0,$ then we have three different alternatives:$\ (i)$ $s$ belongs to some
exceptional set of measure zero, $(ii)$\textsl{\ }$\left(  f^{\ast}\right)
^{^{\prime}}(s)=0,$ or $(iii)$ there is a neighborhood of $s$ such that
$(f^{\ast})^{\prime}(u)$ is not zero, i.e. $f^{\ast}$ is strictly decreasing.
In the two first cases there is nothing to prove. In case alternative $(iii)$
holds then it follows immediately from the properties of the rearrangement
that for a suitable small $h_{0}>0$ we can write
\[
h=\left|  \left\{  f^{\ast}(s+h)<\left|  f\right|  \leq f^{\ast}(s)\right\}
\right|  ,\text{ }0<h<h_{0}.
\]
Therefore for sufficiently small $h$ we can apply Jensen's inequality to
obtain,
\[
\frac{1}{h}\int_{\left\{  f^{\ast}(s+h)<\left|  f\right|  \leq f^{\ast
}(s)\right\}  }\Phi(\left|  \nabla f(x)\right|  )dx\geq\Phi\left(  \frac{1}%
{h}\int_{\left\{  f^{\ast}(s+h)<\left|  f\right|  \leq f^{\ast}(s)\right\}
}\left|  \nabla f(x)\right|  dx\right)  .
\]
Arguing like Talenti \cite{Ta} we thus get
\begin{align*}
\frac{\partial}{\partial s}\int_{\left\{  \left|  f\right|  >f^{\ast
}(s)\right\}  }\Phi(\left|  \nabla f(x)\right|  )dx  &  \geq\Phi\left(
\frac{\partial}{\partial s}\int_{\left\{  \left|  f\right|  >f^{\ast
}(s)\right\}  }\left|  \nabla f(x)\right|  dx\right) \\
&  \geq\Phi\left(  n\gamma_{n}^{1/n}s^{1-1/n}\left(  -f^{\ast}\right)
^{^{\prime}}(s)\right)  ,
\end{align*}
as we wished to show.

\section{\bigskip Proof of the main Theorem \ref{teoA}\label{markateo}}

For the proof we need a slight extension of the following well known fact
(probably due to Hardy and Calder\'{o}n): if $g$ and $h$ are positive and
decreasing and such that
\[
\int_{0}^{t}g(s)ds\preceq\int_{0}^{t}h(s)ds,\forall t>0,
\]
then for any r.i. norm $X$ we have
\[
\left\|  g\right\|  _{X}\preceq\left\|  h\right\|  _{X}.
\]
We extend this result as follows

\begin{lemma}
\label{balboa}Let $f$ and $g$ be two positive functions on the half line.
Moreover, suppose that there exists a real number $\alpha$ such that the
function $t^{\alpha}f(t)$ is monotone (increasing or decreasing). Then, for
any r.i.\ space $X$ with lower Boyd index $\alpha_{X}>0$, there exists a
constant $C=C(\alpha,X)$ such that if $\int_{0}^{t}f(s)\,ds\leq\int_{0}%
^{t}g(s)\,ds,$ holds for all\ $t>0,$ then
\[
\Vert f\Vert_{X}\leq C\Vert g\Vert_{X}.
\]
\end{lemma}

\begin{proof}
Let $Pg(t)=\frac{1}{t}\int_{0}^{t}g(s)ds$ and its adjoint $Qg(t)=\int
_{t}^{\infty}g(s)\frac{ds}{s}$ be the usual Hardy operators (notice that $Q$
is a positive operator and that $Qg(t)$ is a decreasing function). Then,
applying the operator $Q$ to the inequality $Pf(t)\leq Pg(t) $, and using the
fact that $Q\circ P=P\circ Q$, we obtain
\[
\int_{0}^{t}Qf(s)\,ds\leq\int_{0}^{t}Qg(s)\,ds,\text{ \ for all }\ t>0.
\]
Since the integrated functions are decreasing we can apply the usual
Hardy-Calder\'{o}n Lemma (see the discussion preceeding this lemma) to obtain
\[
\Vert Qf\Vert_{X}\leq\Vert Qg\Vert_{X}.
\]
Moreover, since $\alpha_{X}>0,$ we can continue with
\begin{equation}
\Vert Qf\Vert_{X}\leq c_{X}\Vert Q\Vert_{X\rightarrow X}\Vert g\Vert_{X}.
\label{derecha}%
\end{equation}
To estimate the left hand side of (\ref{derecha}) from below we assume first
that the function $t^{\alpha}f(t)$ is increasing. If $\alpha\neq0,$ then
\[
Qf(t)\geq\int_{t}^{2t}s^{\alpha}f(s)s^{-\alpha}\frac{ds}{s}\geq t^{\alpha
}f(t)\int_{t}^{2t}s^{-\alpha-1}ds=\frac{1-2^{-\alpha}}{\alpha}\,f(t).
\]
While if $\alpha=0$ then we readily see that $Qf(t)\geq\frac{1}{2}f(t).$
Similarly, if the function $t^{\alpha}f(t)$ is decreasing, $\alpha\neq0,$
then
\[
Qf(t)\geq\int_{t}^{2t}s^{\alpha}f(s)s^{-\alpha}\frac{ds}{s}\geq(2t)^{\alpha
}f(2t)\int_{t}^{2t}s^{-\alpha-1}ds=\frac{2^{\alpha}-1}{\alpha}\,f(2t).
\]
While if $\alpha=0$ then we readily see that $Qf(t)\geq\frac{1}{2}f(2t).$
Thus, if $t^{\alpha}f(t)$ is monotone, we have
\begin{equation}
\left\|  f\right\|  _{X}\leq C(\alpha)\Vert Qf\Vert_{X}. \label{izquierda}%
\end{equation}
Combining (\ref{izquierda}) and (\ref{derecha}) the desired result follows.
\end{proof}

\bigskip

We may now proceed with the proof of Theorem \ref{teoA}

\begin{proof}
In Section \ref{markao} we have proved the implication $(i)\rightarrow(ii).$

$(ii)\rightarrow(iii)$. Using the Hardy operator $P$ we rewrite (\ref{teoA2})
as
\[
P(s^{-\frac{1}{n}}[f^{\ast\ast}(s)-f^{\ast}(s)])(t)\preceq P(\left|  \nabla
f\right|  ^{\ast}(s))(t).
\]
Let $h(s)=$ $s^{-\frac{1}{n}}[f^{\ast\ast}(s)-f^{\ast}(s)],$ and $g(s)=\left|
\nabla f\right|  ^{\ast}(s),$ and note that $s^{1+1/n}h(s)=s[f^{\ast\ast
}(s)-f^{\ast}(s)]=\int_{f^{\ast}(s)}^{\infty}\lambda_{f}(u)du$ (draw a
picture!) is increasing. Therefore, by lemma \ref{balboa}, we find that
\[
\left\|  s^{-1/n}(f^{\ast\ast}(s)-f^{\ast}(s))\right\|  _{X}\preceq\left\|
\left|  \nabla f\right|  \right\|  _{X},
\]
as we wished to show.

$(iii)\rightarrow(iv)$ Let $X=L^{1},$ then (\ref{final}) reads
\[
\int_{0}^{\infty}s^{1-\frac{1}{n}}[f^{\ast\ast}(s)-f^{\ast}(s)]\frac{ds}%
{s}\preceq\left\|  \left|  \nabla f\right|  \right\|  _{1},
\]
and the result follows since formally\footnote{The fact that the integrated
term vanishes can be easily justified by a familiar limiting argument.}
integrating by parts yields
\begin{align*}
\int_{0}^{\infty}s^{1-\frac{1}{n}}[f^{\ast\ast}(s)-f^{\ast}(s)]\frac{ds}{s}
&  =[1-1/n]\int_{0}^{\infty}f^{\ast\ast}(s)s^{1-1/n}\frac{ds}{s}\\
&  =[1-1/n]\left\|  f\right\|  _{L^{n/(n-1),1}}.
\end{align*}

$(iv)\rightarrow(i)$ This is of course trivial since
\[
W_{0}^{1,1}(\mathbb{R}^{n})\subset L^{n/(n-1),1}(\mathbb{R}^{n})\subset
L^{n/(n-1),\infty}(\mathbb{R}^{n}).
\]
\end{proof}

\section{Sobolev Inequalities in r.i. spaces\label{evg}}

In this section we give a self contained approach to the theory of Sobolev
inequalities in the setting of r.i. spaces. Our results provide optimal
results all the way to the borderline cases.

We recall briefly the basic definitions and conventions we use from the theory
of rearrangement-invariant (r.i.) spaces and refer the reader to \cite{BS} for
a complete treatment.

Let $\Omega$ be a domain in $\mathbb{R}^{n}.$ A Banach function space
$X(\Omega)$ is called a r.i. space if $g\in X(\Omega)$ implies that all
functions $f$ \ with the same decreasing rearrangement, $f^{\ast}=g^{\ast}, $
also belong to $X(\Omega),$ and, moreover, $\Vert f\Vert_{X(\Omega)}=\Vert
g\Vert_{X(\Omega)}$. Let us assume that we define $f(x)=0$ whenever
$x\in\mathbb{R}^{n}\setminus\Omega$, then any r.i.\ space $X(\Omega)$ can be
``reduced'' to one-dimensional space (which by abuse of notation we will still
denote by $X),$ $X=X(0,\left|  \Omega\right|  )$ consisting of all
$g:(0,\left|  \Omega\right|  )\mapsto R$ such that $g^{\ast}(t)=f^{\ast}(t)$
for some function $f\in X(\Omega)$. We shall further assume that our r.i.
spaces satisfy the so-called \textit{Fatou property}, i.e., for any sequence
of functions $f_{k}\rightarrow f$ a.e$,$ $f_{k}\in X,$ and such that $\sup
_{k}\Vert f_{k}\Vert_{X}\leq M$, it follows that $f\in X$ and $\Vert
f\Vert_{X}\leq\liminf\Vert f_{k}\Vert_{X}$.

The upper and lower Boyd indices\footnote{In terms of the Hardy operators
defined by
\[
Pf(t)=\frac{1}{t}\int_{0}^{t}f(s)ds;\text{ \ \ \ }Q_{a}f(t)=\frac{1}{t^{a}%
}\int_{t}^{\infty}s^{a}f(s)\frac{ds}{s},\text{ \ \ }0\leq a<1;
\]
$P$ (resp. $Q_{a}$) is bounded on $X$ if and only if $\beta_{X}<1$ (resp.
$a<\alpha_{X}$) (see for example \cite[Chapter 3]{BS}). Notice that if $a=0,$
$Q_{0}=Q.$} associated with a r.i. space $X$ are defined by
\[
{\beta}_{X}=\inf\limits_{s>1}\dfrac{\ln h_{X}(s)}{\ln s}\text{ \ \ and
\ \ }\alpha_{X}=\sup\limits_{s<1}\dfrac{\ln h_{X}(s)}{\ln s},
\]
where $h_{X}(s)$ denotes the norm of the dilation operator, i.e.
\[
h_{X}(s)=\sup\limits_{f\in X}\dfrac{\left\|  f^{\ast}(\frac{s}{.})\right\|
_{X(0,|\Omega|)}}{\left\|  f^{\ast}\right\|  _{X(0,|\Omega|)}},s>0.
\]

Furthermore we shall assume, essentially without loss, that the spaces we
consider are separable, and unless otherwise specified we shall also assume
that we work on $\mathbb{R}^{n}.$ However, whenever appropriate, we shall
briefly indicate the necessary modifications to treat more general regular domains.

The results and the proofs of this section are similar to those of the papers
\cite{MP} and \cite{P}, however in our present treatment we have no
restrictions on the upper Boyd index $\beta_{X}$.

We record the following elementary result for Hardy operators (cf. \cite{MP}).

\begin{lemma}
\label{cincouno}Let $X$ be a r.i.\ space with the lower Boyd index $\alpha
_{X}>\alpha\geq0$. Then

(i)
\[
\Vert t^{-\alpha}Qf(t)\Vert_{X}\leq C(\alpha,X)\,\Vert t^{-\alpha}%
f(t)\Vert_{X}.
\]

(ii) If $f^{\ast\ast}(\infty)=0,$ then
\[
\Vert t^{-\alpha}f^{\ast\ast}(t)\Vert_{X}\leq C(\alpha,X)\,\Vert t^{-\alpha
}[f^{\ast\ast}(t)-f^{\ast}(t)]\Vert_{X}.
\]
\end{lemma}

\begin{proof}
Both assertions can be found in \cite{MP}. For example see [\cite{MP}, Lemma
2.5] for (i). To prove (ii) use the Fundamental theorem of Calculus to write
$f^{\ast\ast}(t)=\int_{t}^{\infty}(f^{\ast\ast}(s)-f^{\ast}(s))\frac{ds}{s}$
and apply (i).
\end{proof}

We use the notation
\[
\left|  D^{k}f\right|  =\left(  \sum_{\left|  \alpha\right|  =k}\left|
D^{\alpha}f\right|  ^{2}\right)  ^{1/2}.
\]

\begin{theorem}
Let $X$ be a r.i.\ space with $\alpha_{X}>\frac{k-1}{n}$ for some $k\in
N,\,k<n$. Then there exists a constant $C>0,$ such that
\begin{equation}
\Vert t^{-k/n}[f^{\ast\ast}(t)-f^{\ast}(t)]\Vert_{X}\leq C\,\Vert|D^{k}%
f|\Vert_{X},\text{ }f\in C_{0}^{\infty}(\mathbb{R}^{n}). \label{5.2}%
\end{equation}
\end{theorem}

\begin{proof}
When $k=1$ the condition on $\alpha_{X}$ is simply $\alpha_{X}>0,$ therefore
(\ref{5.2}) for $k=1$ was proved in Theorem \ref{teoA} (iii). We prove the
case $k>1$ by induction. Consider first the case $k=2,$ in which case may
assume that $\alpha_{X}>1/n$. Using (\ref{tres}) we get
\[
\Vert t^{-2/n}[f^{\ast\ast}(t)-f^{\ast}(t)]\Vert_{X}\preceq\Vert
t^{-1/n}|\nabla f|^{\ast\ast}(t)\Vert_{X}.
\]
Applying Lemma \ref{cincouno} with $\alpha=1/n$ we can continue with
\[
\Vert t^{-1/n}|\nabla f|^{\ast\ast}(t)\Vert_{X}\preceq\Vert t^{-1/n}[|\nabla
f|^{\ast\ast}(t)-|\nabla f|^{\ast}(t)]\Vert_{X}.
\]
At this point we apply the case $k=1$ to the right hand side to obtain
\[
\Vert t^{-1/n}[|\nabla f|^{\ast\ast}(t)-|\nabla f|^{\ast}(t)]\Vert_{X}%
\preceq\Vert|\nabla\left|  \nabla f\right|  |\Vert_{X}\preceq\Vert
|D^{2}f|\Vert_{X}.
\]
Combining these inequalities thus proves the desired result for the case
$k=2.$ The general case is obtained with the same argument. Indeed, assuming
the inequality is valid for $k-1,$ we can write
\begin{align*}
\Vert t^{-k/n}[f^{\ast\ast}(t)-f^{\ast}(t)]\Vert_{X}  &  \preceq\Vert
t^{-(k-1)/n}|\nabla f|^{\ast\ast}(t)\Vert_{X}\\
&  \leq C_{k-1}\Vert t^{-(k-1)/n}[|\nabla f|^{\ast\ast}(t)-|\nabla f|^{\ast
}(t)]\Vert_{X}\\
&  \leq C_{k}\Vert\,|\nabla|\nabla^{k-1}f||\,\Vert_{X}\\
&  \leq C\,\Vert\,|D^{k}f|\Vert_{X},
\end{align*}
and the result follows.
\end{proof}

To formulate necessary conditions we consider the linear integral
operators\footnote{In the case of domains $\Omega$ one needs to consider
likewise the operators $\tilde{H}_{k/n}g(t)=\int_{t}^{\left|  \Omega\right|
}s^{k/n}g(s)\frac{ds}{s}.$} (cf. \cite{CP}, \cite{EKP}, \cite{MP}, \cite{MM},
\cite{KP1} and the references therein)%
\[
H_{k/n}g(t)=\int_{t}^{\infty}s^{k/n}g(s)\frac{ds}{s}.
\]
The next result was recorded in \cite{P} but with the restriction $\beta
_{X}<1,$ the restriction was later removed in \cite{KP2} but with a rather
complicated proof. \ Our proof provides a considerable simplification.

\begin{theorem}
\label{cincotres}Let $k\in N,\,k<n,$ and let $X$ be a r.i.\ space such that
$\alpha_{X}>\frac{k-1}{n};$ and let $Y$ be another r.i.\ space. Then there
exists a constant $C>0$ such that $\Vert f\Vert_{Y}\leq C\,\Vert
\,|D^{k}f|\,\Vert_{X}$ for all $f\in C_{0}^{\infty}(\mathbb{R}^{n})$ if and
only if $H_{k/n}$ is a bounded operator from $X\rightarrow Y$.
\end{theorem}

\begin{proof}
Suppose that $H_{k/n}$ is a bounded operator, $H_{k/n}:X\rightarrow Y$. Let
$f\in C_{0}^{\infty}(\mathbb{R}^{n}),$ and define $g(t)=t^{-k/n}[f^{\ast\ast
}(t)-f^{\ast}(t)]$, then
\[
H_{k/n}g(t)=\int_{t}^{\infty}[f^{\ast\ast}(s)-f^{\ast}(s)]\frac{ds}%
{s}=Q(f^{\ast\ast}-f^{\ast})(t)=f^{\ast\ast}(t).
\]
Therefore,
\begin{align*}
\Vert f\Vert_{Y}  &  \leq\Vert f^{\ast\ast}\Vert_{Y}=\Vert H_{k/n}g\Vert_{Y}\\
&  \leq\Vert H_{k/n}\Vert_{X\rightarrow Y}\Vert t^{-k/n}[f^{\ast\ast
}(t)-f^{\ast}(t)]\Vert_{X}\\
&  \preceq\Vert|D^{k}f|\Vert_{X}\text{ (by (\ref{5.2})).}%
\end{align*}

To prove the converse we consider first the case $k=1.$ Suppose that $Y$ is a
r.i.\ space such that $\Vert f\Vert_{Y}\leq C\,\Vert\,|\nabla f|\,\Vert_{X}$
for all admissible $f$. Let $g$ be an arbitrary non-negative function from
$X$; we must show that the function $u$ defined by
\[
u(t)=H_{1/n}g(t)=\int_{t}^{\infty}s^{1/n}g(s)\frac{ds}{s},
\]
belongs to $Y.$ Note that $u^{\prime}(t)=\frac{1}{n}t^{1/n-1}g(t),$ therefore
if we define $f(x)=u(t)$ with $t=|x|^{n}$, we see that $|\nabla
f(x)|=nt^{1-1/n}|u^{\prime}(t)|$, whence $\left|  \nabla f(x)\right|  =ng(t).$
It follows that
\begin{align*}
\left\|  u\right\|  _{Y}  &  \simeq\left\|  f\right\|  _{Y}\preceq
\Vert\,|\nabla f|\,\Vert_{X}\text{ (by hypothesis)}\\
&  =Cn\Vert g\Vert_{X},
\end{align*}
as we wished to show. Suppose now that $k>1.$ Repeating the previous argument
$k$ times leads to the conclusion that the operators $(H_{1/n})^{k}$ are
bounded, $(H_{1/n})^{k}:X\rightarrow Y.$ In particular, there exists an
absolute constant $c>0$ such that $\left\|  (H_{1/n})^{k}g\right\|  _{Y}\leq
c\left\|  g\right\|  _{X}.$ To prove that $H_{k/n}$ is a bounded operator
$H_{k/n}:X\rightarrow Y,$ we compare $H_{k/n}$ with $(H_{1/n})^{k}$. By
induction we find
\[
(H_{1/n})^{k}g(t)=n^{k-1}\int_{t}^{\infty}s^{1/n}(s^{1/n}-t^{1/n}%
)^{k-1}g(s)\frac{ds}{s}.
\]
It follows by direct calculations that there exist constants $c_{m},a_{n}$
such that
\begin{equation}
H_{k/n}g(t)=\sum_{m=0}^{k-1}c_{m}\,t^{m/n}(H_{1/n})^{k-m}g(t) \label{5.4}%
\end{equation}%
\[
(H_{1/n})^{k}g(\frac{t}{2})\geq(a_{n})^{m}\,t^{m/n}(H_{1/n})^{k-m}g(t),\qquad
m=1,2,\,\ldots,k-1.
\]

Since the operators $(H_{1/n})^{k}$ are bounded and the dilation operator is
bounded on any r.i.\ space, it follows that
\begin{align*}
\left\|  t^{m/n}(H_{1/n})^{k-m}g(t)\right\|  _{Y}  &  \preceq\left\|
(H_{1/n})^{k}g\right\|  _{Y}\\
&  \preceq\left\|  g\right\|  _{X}.
\end{align*}
Whence from (\ref{5.4}) we obtain that
\[
\left\|  H_{k/n}g\right\|  _{Y}\preceq\left\|  g\right\|  _{X},
\]
as we wished to show.

\begin{remark}
A similar proof of the necessity part is given in \cite{KP2}.
\end{remark}
\end{proof}

\begin{corollary}
Let $k\in N,\,k<n,$ and let $X$ be a r.i.\ space such that $\alpha_{X}%
>\frac{k-1}{n};$ and let $Y$ be another r.i.\ space. Then there exists a
constant $C>0$ such that $\Vert f\Vert_{Y}\leq C\,\Vert\,|D^{k}f|\,\Vert_{X}$
for all $f\in C_{0}^{\infty}(\mathbb{R}^{n})$ if and only if
\[
\Vert f\Vert_{Y}\preceq\Vert t^{-k/n}[f^{\ast\ast}(t)-f^{\ast}(t)]\Vert
_{X},\text{ }f\in C_{0}^{\infty}(\mathbb{R}^{n}).
\]
\end{corollary}

\begin{proof}
Suppose that $\Vert f\Vert_{Y}\leq C\,\Vert\,|D^{k}f|\,\Vert_{X}$ for all
$f\in C_{0}^{\infty}(\mathbb{R}^{n}).$ Let $f\in C_{0}^{\infty}(\mathbb{R}%
^{n}),$ then by (\ref{5.2}) $t^{-k/n}[f^{\ast\ast}(t)-f^{\ast}(t)]\in X$ and
consequently by Theorem \ref{cincotres} we get,
\[
\Vert H_{k/n}(t^{-k/n}[f^{\ast\ast}(t)-f^{\ast}(t)])\Vert_{Y}\preceq\left\|
t^{-k/n}[f^{\ast\ast}(t)-f^{\ast}(t)]\right\|  _{X}.
\]
On the other hand, since
\[
H_{k/n}(t^{-k/n}[f^{\ast\ast}(t)-f^{\ast}(t)])=Q(f^{\ast\ast}-f^{\ast
})=f^{\ast\ast},
\]
we see that
\[
\left\|  f\right\|  _{Y}\leq\left\|  f^{\ast\ast}\right\|  _{Y}\preceq\left\|
t^{-k/n}[f^{\ast\ast}(t)-f^{\ast}(t)]\right\|  _{X}%
\]
as we wished to show.
\end{proof}

The previous discussion provides a method to construct the optimal range space
for a Sobolev inequality. Indeed, let $X$ be a r.i. space with $\alpha
_{X}>\frac{k-1}{n},$ and let the Sobolev space $W_{0}^{k,X}=W_{0}%
^{k,X}(\mathbb{R}^{n})$ be defined to be the closure of $C_{0}^{\infty
}(\mathbb{R}^{n})$ under the norm $\Vert\,|D^{k}f|\,\Vert_{X}.$ Then the
optimal target space $Y$ for the embedding $W_{0}^{k,X}\subset Y$ is given by
the condition
\begin{equation}
\left\|  f\right\|  _{Y}=\Vert t^{-k/n}[f^{\ast\ast}(t)-f^{\ast}(t)]\Vert
_{X}<\infty. \label{facil}%
\end{equation}
However, the space $Y$ defined by (\ref{facil}) may not give a linear function
space. For example, if $X=L^{n/k},k<n,$ then the optimal range space for
Sobolev's inequality is given by the condition (cf. \cite{BMR}, \cite{MP})
\[
\left\|  f\right\|  _{Y}=\left\{  \int_{0}^{\infty}(t^{-k/n}[f^{\ast\ast
}(t)-f^{\ast}(t)])^{n/k}dt\right\}  ^{k/n}=\left\|  f\right\|  _{L(\infty
,n/k)}<\infty,
\]
which is not a linear space. On the other hand, away from the borderline case
(i.e. with a more restrictive condition on the lower Boyd index) it is easy to
see that (\ref{facil}) is equivalent to a r.i. Banach space.

In what follows it will be useful to formally define when a Sobolev embedding
is optimal.

\begin{definition}
Let $X,Y$ be r.i. spaces such that we have a continuous embedding $W_{0}%
^{k,X}\subset Y.$ We shall say that $W_{0}^{k,X}\subset Y$ is optimal if given
any other r.i. $Z$ such that $W_{0}^{k,X}\subset Z,$ it follows that $Y\subset
Z$ continuously.
\end{definition}

\begin{corollary}
\label{canseco}Let $X$ be a r.i.\ space with $\alpha_{X}>\frac{k}{n}$ for some
$k\in N,\,k<n$, and let $Y$ be the r.i.\ space defined by the norm $\Vert
f\Vert_{Y}=\Vert t^{-k/n}f^{\ast\ast}(t)\Vert_{X}$. Then $W_{0}^{k,X}\subset
Y,$ and the embedding is optimal.
\end{corollary}

\begin{proof}
By Lemma \ref{cincouno} with $\alpha=k/n,$%
\[
\Vert t^{-k/n}f^{\ast\ast}\Vert_{X}\preceq\Vert t^{-k/n}[f^{\ast\ast
}(t)-f^{\ast}(t)]\Vert_{X}.
\]
The result now follows from the previous Corollary.
\end{proof}

We conclude discussing how our results can be applied to simplify the study of
compactness of Sobolev embeddings in the setting of r.i. spaces. For the study
of compactness it is natural to restrict oneself to bounded domains $\Omega,$
and henceforth all spaces will be assumed to be based on a bounded domain
$\Omega$ with smooth boundary.

In the study of compactness we will use the following characterization of
compact sets (cf. \cite{P1} and the references therein):

\begin{lemma}
\label{pupu}Let $Z$ be a r.i. space and let $H\subset Z$ be a bounded.set.
Then $H$ is compact in $Z$ iff $H$ is compact in measure and $H$ has
absolutely equicontinuous norm\footnote{Recall that a set $H\subset Z$ is
absolutely equicontinuous in norm if $\forall\varepsilon>0$ $\exists\delta>0$
such that if $\left|  D\right|  <\delta$ then $\left\|  f\chi_{D}\right\|
_{Z}<\varepsilon.$}.
\end{lemma}

In order to use the results of this paper we recall the connection between
optimal embeddings and compactness. Indeed, it is known from the classical
$L^{p}$ theory that optimal Sobolev embeddings are not compact. Pustylnik
\cite{P1}, has recently extended this result and, most importantly for our
purposes, quantified the lack of compactness of optimal embeddings. More
precisely, we have the following (cf. \cite{P1})

\begin{lemma}
\label{pu}Suppose that $W_{0}^{k,X}\subset Y$ is optimal, and let $Z$ be a
r.i. space such that $W_{0}^{k,X}\subset Z$ is compact. Then the inclusion
$Y\subset Z$ is absolutely continuous\footnote{This means that every bounded
set $H\subset Y$ is absolutely continuous in norm in $Z.$}.
\end{lemma}

We also note for future use that by an easy case of the Rellich-Kondrachov
theorem, the embedding $W_{0}^{1,L^{1}}\subset L^{1}$ is compact. Therefore,
since for any r.i. space $X$ we have $W_{0}^{k,X}\subset W_{0}^{1,L^{1}%
}\subset L^{1},$ we see that all bounded sets in $W_{0}^{k,X}$ are compact in
measure. Consequently to verify that an embedding $W_{0}^{k,X}\subset Z$ is
compact it is only necessary to verify that bounded sets in $W_{0}^{k,X}$ have
absolutely continuous norm in $Z.$

With these preliminaries at hand we shall now provide our proof of the
compactness result recently obtained in \cite{P1} and \cite{KP3} with
different but long and complicated methods of proof.

\begin{theorem}
\label{marke}Let $X,Z$ be r.i. spaces with $\alpha_{X}>\frac{k-1}{n}$ and such
that $W_{0}^{k,X}\subset Z.$ Then the embedding $W_{0}^{k,X}\subset Z$ is
compact if and only if $\tilde{H}_{\frac{n}{k}}$ is a compact operator
$\tilde{H}_{\frac{n}{k}}:X\rightarrow Z,$ here $\tilde{H}_{\frac{n}{k}%
}f(t)=:\int_{t}^{\left|  \Omega\right|  }s^{k/n}f(s)\frac{ds}{s}$.
\end{theorem}

\begin{proof}
Suppose first that the embedding $W_{0}^{k,X}\subset Z$ is compact and
consider the optimal embedding $W_{0}^{k,X}\subset Y$ provided by
(\ref{facil}) or by Corollary \ref{canseco}. It follows readily, by a suitable
modified version of Theorem \ref{cincotres} for bounded domains, that
$\tilde{H}_{\frac{n}{k}}:X\rightarrow$ $Y$ is bounded. It is easy to see that
this implies that $\tilde{H}_{\frac{n}{k}}$ sends bounded sets $A\subset X$
into sets $\tilde{H}_{\frac{n}{k}}(A)$ which are compact in measure. Moreover,
by Pustylnik's Lemma \ref{pu} , the embedding $Y\subset Z$ is absolutely
equicontinuous and since we obviously can factor $\tilde{H}_{\frac{n}{k}%
}:X\rightarrow$ $Y\subset Z,$ we see that $\tilde{H}_{\frac{n}{k}%
}:X\rightarrow Z$ also maps bounded sets into sets that are absolutely
equicontinuous. Therefore, from the compactness criteria given by Lemma
\ref{pupu}, we find that $\tilde{H}_{\frac{n}{k}}:X\rightarrow Z$ is a compact operator.

Conversely, suppose that $\tilde{H}_{\frac{n}{k}}:X\rightarrow Z$ is a compact
operator, and let $A$ be a bounded set in $W_{0}^{k,X}.$ By the definition of
$W_{0}^{k,X}$ we may assume without loss that $A\subset C_{0}^{\infty}.$ As
pointed out above $A$ is automatically compact in measure, therefore, by Lemma
\ref{pupu}, to prove that $A$ is compact in $Z$ it remains to verify that $A$
has absolutely equicontinuous norm. Define $\tilde{A}=\{\tilde{f}:\tilde
{f}(t)=t^{-k/n}[f^{\ast\ast}(t)-f^{\ast}(t)],$ $f\in A\}.$ By (\ref{5.2}),
$\tilde{A}$ is a bounded set in $X=X(0,\left|  \Omega\right|  ),$ therefore
$\tilde{H}_{\frac{n}{k}}(\tilde{A})$ is compact in $Z,$ in particular it has
absolutely equicontinuous norm, $\lim_{a\rightarrow0}\sup_{f\in\tilde{A}%
}\left\|  \tilde{H}_{\frac{n}{k}}\tilde{f}\chi_{(0,a)}\right\|  _{Z}=0.$
Moreover, since
\[
\tilde{H}_{\frac{n}{k}}\tilde{f}\geq f^{\ast\ast}\geq f^{\ast}%
\]
it follows that
\[
\lim_{a\rightarrow0}\sup_{f\in A}\left\|  f\chi_{(0,a)}\right\|  _{Z}=0,
\]
and consequently $A$ has absolutely equicontinuous norm as we wished to show.
\end{proof}


\begin{thebibliography}{99}
\bibitem{BaCLS}D. Bakry, T.\ Coulhon, M. Ledoux, and L. Saloff-Coste, Sobolev
inequalities in disguise, Indiana Univ. Math. J. 44 (1995), 1033-1074.

\bibitem{BMR}J. Bastero, M. Milman and F. Ruiz, \textsl{A note on }%
$L(\infty,q)$\textsl{\ spaces and Sobolev embeddings}, Indiana Univ. Math. J.
52 (2003), 1215-1230.

\bibitem{BS}C. Bennett and R. Sharpley, \textsl{Interpolation of Operators},
Academic Press, Boston\textbf{, }1988

\bibitem{CP}M. Cwikel and E. Pustylnik, \textsl{Sobolev type embeddings in the
limiting case}, J. Fourier Anal. Appl. 4 (1998), 433-446.

\bibitem{EKP}D. E. Edmunds, R. Kerman and L. Pick, \textsl{Optimal Sobolev
embeddings involving rearrangement invariant quasi-norms,} J. Funct. Anal. 170
(2000), 307-355.

\bibitem{Ha}P. Hajlasz,\textsl{\ Sobolev inequalities, truncation method, and
John domains}, Papers in Analysis, Rep. Univ. Jyv\"{a}skyl\"{a} Dep. Math.
Stat. 83, Univ. Jyv\"{a}skyl\"{a}, Jyv\"{a}skyl\"{a}, 2001, pp 109-126.

\bibitem{KM}J. Kalis and M. Milman,\textsl{\ Symmetrization and sharp Sobolev
inequalities in metric spaces}, preprint.

\bibitem{KP1}R. Kerman and L. Pick, \textsl{Optimal Sobolev imbedding spaces},
Preprint MATH-KMA-2005/161, Charles Univ., Prague 2005, pp. 1-19.

\bibitem{KP2}R. Kerman and L. Pick, \textsl{Optimal Sobolev imbeddings}, Forum
Math. 18 (2006), 535-579.

\bibitem{KP3}R. Kerman and L. Pick, \textsl{Compactness of Sobolev imbeddings
involving rearrangement-invariant spaces}, Preprint 2006, pp. 1-26.

\bibitem{Ko}V. I. Kolyada, \textsl{Rearrangements of functions and embedding
theorems}, Uspekhi Mat. Nauk 44 (1989), 61-95 [English transl. in Russian
Math. Surveys 44 (1989), 73-117.

\bibitem{LZ}W. A. J. Luxemburg and A. C. Zaanen, \textsl{Compactness of
integral operators on Banach function spaces}, Math.\ Ann. 149 (1963), 150-180.

\bibitem{MM}J. Mart\'{i}n and M. Milman, \textsl{Higher order symmetrization
inequalities and applications}, J. Math. Anal. Appl., to appear.

\bibitem{MM1}J. Martin and M. Milman, \textsl{Symmetrization inequalites and
Sobolev embeddings}, Proc. Amer. Math. Soc. \textbf{134} (2006), 2335-2347.

\bibitem{Ma}V. G. Maz'ya, \textsl{Sobolev Spaces, }Springer-Verlag, New York, 1985.

\bibitem{MP}M. Milman and E. Pustylnik, \textsl{On sharp higher order Sobolev
embeddings}, Comm. Contemp. Math. 6 (2004), 495-511.

\bibitem{Po}S. Poornima, An embedding theorem for the Sobolev spaces
$W^{1,1}(\mathbb{R}^{n}),$ Bull. Sci. Math. 107 (1983), 253-259.

\bibitem{P}E. Pustylnik, \textsl{Sobolev type inequalities in ultrasymmetric
spaces with applications to Orlicz-Sobolev embeddings}, Funct.\ Spaces Appl. 3
(2005), 183-208.

\bibitem{P1}E. Pustylnik, \textsl{On compactness of Sobolev embeddings}, Forum
Math. 18 (2006), 839--852.

\bibitem{St}E. Stein, \textsl{Singular integrals and differentiability
properties of functions}, Princeton University Press, New Jersey, 1970.

\bibitem{Ta}G. Talenti,\textsl{\ Inequalities in rearrangement-invariant
function spaces}, Nonlinear Analysis, Function Spaces and Applications,
Prometheus, Prague vol. 5, 1995, pp. 177-230.

\bibitem{Tar}L. Tartar, Imbedding theorems of Sobolev spaces into Lorentz
spaces, Boll. Unione Mat. Ital. Sez B Artic. Ric. Mat. (8) 1 (1998), 479-500.
\end{thebibliography}
\end{document}